\documentclass[a4paper,12pt]{amsart}
\usepackage{amsmath,amsthm,amsbsy, amstext, amscd, amsxtra,amsopn, latexsym,amssymb,graphicx,graphics}

\def\shorttitle{vanishing theorem}
\def\shortauthor{{\it  }}
\newtheorem{theorem}{Theorem}[section]{\bf}{\it}
\newtheorem{proposition}{Proposition}[section]{\bf}{\it}
\newtheorem{lemma}{Lemma}[section]{\bf}{\it}
\newtheorem{remark}{Remark}[section]{\it}{\rm}
\newtheorem*{Proof}{Proof}{\it}{\rm}

\renewcommand{\sec}[1]{\setcounter{equation}{0}\section{#1}}

\newcommand{\bslide}{\begin{slide}}
\newcommand{\eslide}{\end{slide}}
\newcommand{\frameb}{\begin{frame}}
\newcommand{\framee}{\end{frame}}
\newcommand{\bitemi}{\begin{itemize}}
\newcommand{\eitemi}{\end{itemize}}
\newcommand{\bcenter}{\begin{center}}
\newcommand{\ecenter}{\end{center}}
\newcommand{\bblock}{\begin{block}}
\newcommand{\eblock}{\end{block}}
\newcommand{\benum}{\begin{enumerate}}
\newcommand{\eenum}{\end{enumerate}}

\newcommand{\riem}{Riemannian }

\newcommand{\wrt}{with respect to }

\newcommand{\s}{\ }

\newcommand{\bthmm}{\begin{thm}}
\newcommand{\ethmm}{\end{thm}}
\newcommand{\blemm}{\begin{lem}}
\newcommand{\elemm}{\end{lem}}
\newcommand{\bppm}{\begin{prop}}
\newcommand{\eppm}{\end{prop}}
\newcommand{\bprfm}{\begin{Proof}}
\newcommand{\eprfm}{\end{Proof}}
\newcommand{\bcorm}{\begin{cor}}
\newcommand{\ecorm}{\end{cor}}

\newcommand{\bthm}{\begin{theorem}}
\newcommand{\ethm}{\end{theorem}}
\newcommand{\blem}{\begin{lemma}}
\newcommand{\elem}{\end{lemma}}
\newcommand{\bpp}{\begin{proposition}}
\newcommand{\epp}{\end{proposition}}
\newcommand{\bprf}{\begin{proof}}
\newcommand{\eprf}{\end{proof}}
\newcommand{\brem}{\begin{remark}}
\newcommand{\erem}{\end{remark}}

\newcommand{\fa}{\forall}

\newcommand{\De}{\Delta}
\newcommand{\de}{\delta}
\newcommand{\na}{\nabla}

\newcommand{\q}{\quad}
\newcommand{\qq}{\qquad}

\newcommand{\pa}{\partial}

\newcommand{\fr}{\frac}

\newcommand{\bsl}{\backslash}

\newcommand{\wg}{\wedge}
\newcommand{\la}{\lambda}

\newcommand{\tm}{\times}

\newcommand{\vep}{\varepsilon}

\newcommand{\gm}{\gamma}

\newcommand{\al}{\alpha}

\newcommand{\tta}{\theta}

\newcommand{\og}{\omega}

\newcommand{\Si}{\Sigma}

\newcommand{\ot}{\otimes}
\newcommand{\op}{\oplus}

\newcommand{\beq}{\begin{eqnarray}}
\newcommand{\eeq}{\end{eqnarray}}
\newcommand{\beqs}{\begin{eqnarray*}}
\newcommand{\eeqs}{\end{eqnarray*}}
\newcommand{\bal}{\begin{align*}}
\newcommand{\eal}{\end{align*}}
\newcommand{\bale}{\begin{aligned}}
\newcommand{\eale}{\end{aligned}}
\newcommand{\bequ}{\begin{equation}}
\newcommand{\eequ}{\end{equation}}
\newcommand{\bequs}{\begin{equation*}}
\newcommand{\eequs}{\end{equation*}}

\newcommand{\bc}{\begin{center}}
\newcommand{\ec}{\end{center}}

\newcommand{\bcase}{\begin{cases}}
\newcommand{\ecase}{\end{cases}}
\newcommand{\bmat}{\begin{matrix}}
\newcommand{\emat}{\end{matrix}}
\newcommand{\bbm}{\begin{bmatrix}}
\newcommand{\ebm}{\end{bmatrix}}
\newcommand{\bpm}{\begin{pmatrix}}
\newcommand{\epm}{\end{pmatrix}}
\newcommand{\bvm}{\begin{vmatrix}}
\newcommand{\evm}{\end{vmatrix}}

\newcommand{\mc}{\mathcal}
\newcommand{\mbbr}{\mathbb{R}}
\newcommand{\mbbc}{\mathbb{C}}

\newcommand{\mbbh}{\mathbb{H}}

\newcommand{\mcr}{\mathcal{R}}

\newcommand{\mfa}{\mathfrak{a}}
\newcommand{\mfb}{\mathfrak{b}}
\newcommand{\mfc}{\mathfrak{c}}
\newcommand{\mfd}{\mathfrak{d}}

\newcommand{\mfg}{\mathfrak{g}}
\newcommand{\mfh}{\mathfrak{h}}

\newcommand{\mfk}{\mathfrak{k}}

\newcommand{\mfp}{\mathfrak{p}}

\begin{document}

\title
[vanishing theorem for irreducible symmetric spaces]
{vanishing theorem for irreducible symmetric spaces of noncompact type}

\author
[Xusheng Liu]{Xusheng Liu}
\address
{School of Mathematical Sciences, Fudan University,
Shanghai 200433, China
}
\email{xshliu@fudan.edu.cn}

\maketitle

\begin{abstract}
We prove the following vanishing theorem.
Let $M$ be an irreducible symmetric space of noncompact type
whose dimension exceeds 2 and  $M\ne SO_0(2,2)/SO(2)\tm SO(2).$
Let $\pi: E\to M$ be any vector bundle,
Then any $E-$valued $L^2$ harmonic 1-form over $M$ vanishes.
In particular we get 
the vanishing theorem for harmonic maps from irreducible symmetric spaces
of noncompact type.
\end{abstract}

\sec{introduction}

It was conjectured by J. H. Sampson \cite{SAM}  that any harmonic map with finite
energy from a complete simply connected Riemannian manifold  with
negative sectional curvature whose dimension exceeds 2 must be
constant. This is valid  for space forms,  but unsolved in
general case. For Cartan-Hadamard manifolds,
In \cite{XIN1} Xin proved a general vanishing theorem as follows. \par
\bthm
\label{thma} Let M be an m-dimensional Cartan-Hadamard manifold
with the sectional curvature  $-a^2 \le K\le 0$ and the Ricci
curvature bounded from above by $-b^2$.  Let f be a harmonic map
from M into any Riemannian manifold with the moderate divergent
energy. If $b\ge 2a$,   then f has to be constant.
\ethm
\par

For Hermitian symmetric spaces of noncompact type, In \cite{XIN2}
Xin proved the following results. \bthm \label{thmb} A harmonic map
of a finite energy from a classical bounded symmetric domain except
$D_{IV}(2)=SO_0(2,2)/SO(2)\tm SO(2)(\cong \mbbh\tm \mbbh)$ to any
Riemannian manifold has to be constant. \ethm

Now in \cite{LIU} by calculating the lower bounds of sectional
curvature of all irreducible symmetric spaces of noncompact type, we
get the following theorem by using Theorem\ref{thma}.

\bthm\label{thmc}
Let $M$ be one of the  irreducible symmetric spaces of
noncompact type in the following cases,
\beqs
&& SL(n,  \mbbr)/SO(n),  n\ge 4;\\
&& SU^*(2n)/Sp(n); \\
&& SU(p, q)/S(U_p\tm U_q),  p+q\ge 4; \\
&&SO_o(p, q)/SO(p)\tm SO(q),  \mbox{for } r=1, p+q\ge 4,\\
&& \hspace{1cm} \mbox{for } r>1, p+q\ge 6, \q \mbox{here } r=min(p, q); \\
&& SO^*(2n)/U(n),  n\ge 3;\\
&& Sp(n, \mbbr)/U(n),  n\ge 3; \\
&&Sp(p, q)/Sp(p)\tm Sp(q); \\
&&  EI, EII, EIII, EIV,  EV, EVI, EVII, EVIII,  EIX,  FI,FII \s \mbox{and } G.
\eeqs
Then any $L^2$ harmonic 1-form vanishes.
\ethm

In this article we prove the following theorem.
\bthm\label{thmm}
Let $M$ be an irreducible symmetric space of noncompact type
whose dimension exceeds 2 and  $M\ne SO_0(2,2)/SO(2)\tm SO(2).$
Let $\pi: E\to M$ be any vector bundle,
Then any $E-$valued $L^2$ harmonic 1-form over $M$ vanishes.
\ethm

\brem
Let $f: M\to N$ be a harmonic map, then $df\in \Lambda^1 f^*TN$ 
is a harmonic 1-form over $M$, from Theorem\ref{thmm} we get
the vanishing theorem for harmonic maps from irreducible symmetric spaces
of noncompact type which generalize Theorem\ref{thma} and \ref{thmb}.
\erem

In section 2 we prove a vanishing theorem for harmonic forms.
In section 3 we give the vanishing theorem for Riemannian symmetric
spaces of noncompact type. In section 4 we calculate
the Hessian of distance function in the cases of
Riemannian symmetric spaces of noncompact type and
give the proof of Theorem\ref{thmm}.

 \sec{the vanishing theorem of  harmonic forms}
Let $M$ be an n dimensional complete \riem manifold,
let $\pi: E\to M$ be a real vector bundle with rank $r$, we
denote $\na$ the metric connection on $E$ which gives rise to
the Levi-Civita connection when restricted to $M$.
\par
Let $\{e_i, 1\le i\le n\}$ be an local orthogonal frame on $M$
with dual coframe fields $\{\tta^i\}.$
We identify the tangent vector field $X$ with a 1-form $\hat X$ via
Riemannian inner product by
\beqs
  \hat X(Y)=<Y,X>.
\eeqs
Then $\hat e_i=\tta^i.$  In the following calculations, we
take the normal coordinate at the given point $x\in M$,
i.e., $\na(e_i)(x)=0.$
\par

Let $\wg^p(M,E)$ be the vector space of all $E$-valued $p$-forms on
M, for $\og\in \wg^p(M,E)$, we have exterior multiplication and
interior product operator \beqs
 &&\vep_X(\og)=\vep_{\hat X}(\og)=\hat X\wg \og, \\
 &&i_X(\og)=i_{\hat X}(\og),\q i_{X}(\og)(X_1,\cdots,X_{p-1})=\og(X,X_1,\cdots,X_{p-1}).
\eeqs
We have the following commutation rules
\beqs
 &&\vep_X\vep_Y+\vep_Y\vep_X=0.\\
 && i_Xi_Y +i_Yi_X=0.\\
 && \vep_Xi_Y+i_Y\vep_X=<X,Y>1, \mbox{ where 1 is the identity operator}.\\
 && \na_X\vep_Y=\vep_Y\na_X+\vep_{\na_XY}.\\
 && \na_Xi_Y=i_Y\na_X+i_{\na_XY}
 \eeqs
 \par

  We can define the adjoint operator  \wrt the inner product on  $\wg^*(M,E)$, then
 \beqs
&& \vep_k^*=\vep_{e_k}^*=i_k=i_{e_k}\\
&& \na_k^*=\na_{e_k}^*=-\na_k
 \eeqs
Now we have the exterior differential operator
 \beqs
 d: \wg^p(M,E)\to \wg^{p+1}(M,E), \q d=\vep_k\na_k=\vep_{\tta^k}\na_{e_k}.
 \eeqs
its  adjoint operator is
 \beqs
 \de=d^*=-i_k\na_k.
 \eeqs
 \par
Let $X$ be a tangent vector field, we have the Lie derivative \wrt $X$
 \beqs
 L_X=d\circ i_X+i_X\circ d.
 \eeqs
We have
 \beqs
 d\circ i_X&=&\vep_k\na_ki_X=\vep_k(i_X\na_k+i_{\na_kX}) \\
 &=&(-i_X\vep_k+<X,e_k>)\na_k+\vep_ki_{\na_kX}\\
 &=&-i_X\vep_k\na_k+\na_X+\vep_k i_{\na_kX}\\
 i_X\circ d&=&i_X\vep_k\na_k \\
 L_X&=&\na_X+\vep_ki_{\na_kX}.
 \eeqs
 \par
In particular, let $f(x)$ be a function on $M$,
we denote  $X=\na f(x)$  its  gradient vector field.
The Hessian of $f$ is
 \beqs
Hess(f)=h_{kl}\tta^k\ot \tta^l=h_{kl}e_k\ot e_l.
\eeqs
Then
\beqs
 \na_k\na f=h_{kl}e_l, L_{\na f}=\na_{\na f}+h_{kl}\vep_ki_l.
 \eeqs
It follows that
 \beqs\label{eq00}
 <L_{\na f}\og,\og>&=&<\na_{\na f}\og, \og>+h_{kl}<\vep_ki_l\og,\og>\\
  &=& \fr 12\na_{\na f}|\og|^2+h_{kl}<i_k\og, i_l\og>\\
  &=& \fr 12<\na f, \na |\og|^2>+H_f(\og,\og).
 \eeqs
where
 \beqs
 H_f(\og,\og)=\sum_{kl}h_{kl}<i_k\og, i_l\og>=\sum <h_{kl}\vep_ki_l\og,\og>.
 \eeqs
 \par
We recall the Green formula,
let $D\subset M$ be a relatively compact domain with smooth boundary,
we denote $\nu$ the unit outward normal vector field on boundary.
We have
  \beqs
 &&\int_D[<d\og,\psi>-<\og, \de\psi>]\\
 &&=\int_{\pa D}<\og, i_\nu \psi>, \q \og\in \wg^p(M,E), \psi\in \wg^{p+1}(M,E).
 \eeqs
We calculate the integral in two ways
 \beqs\label{eq01}
 \int_D <L_{\na f}\og,\og>&=&\int_D <d\circ i_{\na f} \og+i_{\na f}\circ d \og,\og>\\
 &=& \int_D <i_{\na f}\og,\delta \og>+ \int_{\pa_D}<i_{\na f}\og, i_\nu\og>\\
 &&+\int_D<d\og, \vep_{\na f}\og >\\
 \int_D <\na f, \na |\og|^2>&=& \int_D div(|\og|^2 \na f)-\De f|\og|^2\\
 &=& -\int_D \De f |\og|^2+\int_{\pa D}\na_\nu f |\og|^2.
 \eeqs
We get
 \beqs
 &&\int_D [-\De f \fr {|\og|^2}2+H_f(\og,\og) -<i_{\na f}\og,\delta \og> -<d\og, \vep_{\na f}\og >]\\
 &&=\int_{\pa D}[<i_{\na f}\og, i_\nu\og>-\fr{|\og|^2}2\na_\nu f]
 \eeqs
in other words,
 \beqs \label{integral1}
&& \int_D [\De f |\og|^2-2H_f(\og,\og)+2<i_{\na f}\og,\delta \og> +2<d\og, \vep_{\na f}\og >]\\
&&=\int_{\pa D}[-2<i_{\na f}\og, i_\nu\og>+|\og|^2\na_\nu f]
 \eeqs
 \par

Let $o\in M$ be a fixed point and $r(x)=dist(o,x)$ the distance
function from $o$. We denote $B_r=B(o,r)=\{x\in M| r(x)< r\},
S_r=\pa B_r=\{x\in M| r(x) =r\}.$ The Hessian of $r(x)$ is the
second fundamental form \wrt $S_r,$
$$
Hess(r)(X,Y)=<\nabla_X \fr{\pa}{\pa r},Y>,\q  X,Y\in TS_r(x).
$$
Let $\la_1(x)\ge \la_2(x))\cdots \ge  \la_{n-1}(x)$ be the $n-1$ eigenvalues  of $Hess(r)$,
we get the Laplacian of $r$, i.e.,
$$
\Delta r=tr Hess(r)=\sum_{1\le i\le n-1}\la_i(x).
$$

\bthm\label{vanish1}
Let $M$ be a noncompact complete \riem manifold, $0\le p\le n$ be an integer,
let $\pi: E\to M$ be a vector bundle.
Let $o\in M$ be a fixed point with the distance function $r(x)=dist(o,x).$
We order the   eigenvalues  of $Hess(r)$ with multiplicities in the way that
$\la_1(x)\ge \la_2(x))\cdots \ge  \la_{n-1}(x).$
If there exists $R_0\ge 0$, we have for $r(x)\ge R_0$,
\beq \label{conda}
 \sum_{i=1}^p\la_i(x) \le \sum_{i=p+1}^{n-1}\la_i(x).
\eeq
Then any $E$-valued $ L^2$ p-harmonic form on $M$ vanishes.
\ethm

\begin{proof}
Let $\omega\in \wg^p(M,E)$, let $f(x)=r(x), \nu=\na r$.
Since $d\og=\de\og=0$,
From (\ref{integral1}) we get
\beqs
 \int_D [\De r |\og|^2-2H_r(\og,\og)]&=&\int_{\pa D}[-2<i_{\nu}\og, i_\nu\og>+|\og|^2\na_\nu r]\\
 &\le & \int_{\pa D} |\og|^2
 \eeqs
 \par
At $x$ we choose normal frame $\{e_i\}$ with
$e_n=\nu=\na r, e_i\in TS_r, 1\le i\le n-1,$
we suppose $Hess(r)$ is diagonal \wrt $\{e_i\}$, then
 \beqs
 Hess(r)&=&\sum_{1\le i\le n-1}\la_i(x)e_i\ot e_i.\\
 \De r &=&\sum_{1\le i\le n-1}\la_i(x).\\
 Hess_r(\og,\og)&\le &\sum_{1\le i\le p}\la_i <\og,\og>.
 \eeqs
We get
 \beqs
&&\De r|\og|^2-Hess_r(\og,\og)\ge [\fr 12 \Delta r-\sum_{1\le i\le p}\la_i(x)]|\og|^2\\
&=& \fr 12 [\sum_{p+1\le i\le n-1}\la_i(x)-\sum_{1\le i\le p}\la_i(x)]|\og|^2
\eeqs
\par
Under the condition (\ref{conda}), we take $D=B_{R_2}\bsl B_{R_1},
R_0\le R_1<R_2,$ then \beqs 0\le
\int_{S_{R_2}}|\og|^2-\int_{S_{R_1}}|\og|^2. \eeqs We see that
$\int_{S_R}|\og|^2$ is non decreasing \wrt $R$ for $R\ge R_0$. The
unique continuation theorem of harmonic forms says that the harmonic
form $\og$ must vanishes anywhere when $\og$ vanishes in some open
subset of $M$. We see that if $\og\ne 0$, then there exists $R\ge
R_0$ such that $\int_{S_{R}}\fr 12|\og|^2\ge c>0$, which gives rise
to \beqs
 \fr 12 \int_M|\og|^2\ge\int_{R}^{\infty}(\int_{S_r}\fr 12|\og|^2)dr\ge \int_{R}^{\infty} cdr=\infty.
 \eeqs
this contradicts  to the condition of  finite energy.
 \end{proof}

 \brem
In fact we can get the vanishing theorem if we have the estimate
 $
 \De r|\og|^2-Hess_r(\og,\og)\ge 0 $
 outside some compact subset of $M$,
 We can give the similar proof by using the conservation law of stress-energy
 in \cite{JX}.
 The formula (\ref{integral1}) is given at \cite{EF}.
\erem

\sec{the  vanishing theorem in the case of \riem symmetric spaces}

Let M be a simply connected \riem symmetric space of non compact type with dimension $n$.
Let $o\in M$ be a fixed point with the distance function $r(x)=dist(o,x)$.
$\gamma(t): [0,r(x)]\to M $ be the unique geodesic  from o to x. The curvature transformation along
$\gm$ is the self adjoint operator
\beqs
&\mc{R}_{\dot\gm(t)}V=R(V,\dot\gm(t))\dot\gm(t), \\
&<\mc{R}_{\dot\gm(t)}V,W>=R(\dot\gm(t),V,\dot\gm(t),W), \fa V,W\perp \dot\gm(t).
\eeqs
\par
We choose the orthogonal frame $\{e_i,1\le i\le n\}$ along $\gm$
so that $e_n(t)=\dot\gm(t)$.
 Let $-\la_i^2(t), 1\le i\le n-1\ge 0$ be the eigenvalues of $\mc{R}_{\dot\gm(t)}$,
 we order $\la_i(t)$ with multiplicities in the form of
 $\la_1(t)\ge \la_2(t)\ge\cdots \ge \la_{n-1}(t)$.
 Since the curvature tensor of M is parallel, all $\la_i(t)$ are constant along
 $\gm, \la_i(t)\equiv\la_i(0)=\la_i$.
 At x, we have
 \beqs
 \Delta r(x)=\sum_{i=1}^{n-1}\la_i \coth(\la_ir)\\
 Hess(r)=\sum_{i=1}^{n-1}\la_i\coth(\la_i)e_i(r)\ot e_i(r)\\
 Ric(\fr{\pa}{\pa r},\fr{\pa}{\pa r})=Ric(e_n,e_n)=-\sum_{i=1}^{n-1}\la_i^2
 \eeqs
  \begin{theorem}\label{vanish2}
 Let $M$ be a simply connected \riem symmetric space of non compact type,
 let $\pi: E\to M$ be a vector bundle.
let  $0\le p \le n$ be  an integer.
If at some point $o\in M$(hence  at all points),
the eigenvalues  of  curvature transformation operator $\mc{R}_v, v\in T_oM$ be any unit vector at $o$,
 which ordered by $\la_1\ge\la_2\ge\cdots\ge\la_{n-1}$ satisfy
 \beq \label{cond1}
 \sum_{i=1}^p\la_i\le \sum_{i=p+1}^{n-1}\la_i.
 \eeq
Then any $L^2$ $E$-valued  harmonic p-forms vanish.
  In particular, let the Ric curvature be $Ric=-B$, the curvature lower bounds of M is $-A$,
 if
 \beq \label{cond2}
 p(p+1)\le \fr BA,
 \eeq
 then  any $L^2$ harmonic p-form vanishes.
  \end{theorem}
\begin{proof}
\par

If (\ref{cond1}) holds, as $\phi(\la)=\coth(\la r)$ is decreasing \wrt $\la$,
we have
\beqs
\sum_{p+1\le i\le n-1}\la_i\coth(\la_ir)-\sum_{1\le i\le p}\la_i\coth(\la_ir)\ge 0
\eeqs
\par

Then  if (\ref{cond1}) holds, from Theorem \ref{vanish1},
we get any $L^2$ harmonic form vanish.\par
If (\ref{cond2}) holds,
we have
\beqs
\sum_{1\le i\le n-1}\la_i^2\ge p(p+1)\la_1^2\ge (p+1)\sum_{1\le i\le p}\la_i^2.
\eeqs
i.e.,
\beqs
\sum_{p+1\le i\le n-1}\la_i^2\ge p\sum_{1\le i\le p}\la_i^2\ge (\sum_{1\le i\le p} \la_i)^2.
\eeqs
We get
\beqs
(\sum_{p+1\le i\le n-1}\la_i)^2\ge
\sum_{p+1\le i\le n-1}\la_i^2\ge (\sum_{1\le i\le p} \la_i)^2.
\eeqs
Then (\ref{cond1}) holds.
\end{proof}

Remark: for p=1 we recover the theorem 2.2 in \cite{JX}.

\sec{the curvature tensor of irreducible \riem symmetric spaces of
non compact type} Let $M=G/K$ be a \riem symmetric space of non
compact type with Lie algebra decomposition $\mfg=\mfk+\mfp$ and
Carton involution $\theta$. Let $\mfa\subset \mfp$ be a maximal
abelian subspace over $\mbbr$, $rank(M)=dim(\mfa)=r$, we extend
$\mfa$ to a Cartan subalgebra
 of $\mfg$, say $\mfh=\mfh_{\mfk}+\mfa, \mfh_{\mfk}\subset \mfk$.
  Let $\Delta=\Delta(\mfg^{\mbbc}, \mfh^{\mbbc})$ be the corresponding
 root system of complex semisimple Lie algebra $\mfg^{\mbbc}$.
 Let $\Sigma$ be the restricted root system, i.e, consist of restriction  of root of $\Delta$ to $\mfa$.
 For $\la\in \Sigma$, we denote $\mfg_{\la}$ the root space with multiplicity
 $m_{\la}=dim \mfg_{\la}$.
We have the Iwasawa decomposition  $\mfg=\mfk\op\mfa\op\sum_{\la\in\Si^+}\mfg_{\la}$.
Any vector in $\mfp$ is conjugate by $Int(\mfk)$ to a vector in $\mfa$.
Let (X,Y)=B(X,Y) be the killing form of $\mfg$ which induces a left invariant inner product
on M by $<X,Y>=-B(X,\tta Y)$. For any root $\al$, We embed it into $\mfh$ by
$\al(h)=(\al,h), h\in \mfh$. The curvature tensor of M is
\beqs
&&R(X,Y)=[\na_X,\na_Y]-\na_{[X,Y]}=-ad[X,Y]\\
&&R(X,Y,Z,W)=-<R(X,Y)Z,W>\\
&&=([X,Y],[Z,W])=-<[X,Y],[Z,W]>.
\eeqs
\par
Let $h\in \mfa, |h|=1, X\perp h, X\in \mfp$, then we have \beqs
\mcr_hX=R(X,h)h=-[[X,h],h]=-ad^2h(X). \eeqs\par From this, we get
the eigenvalues and the eigenvectors as follows, \beqs
&&(i) 0, \mbox{the orthogonal complement of } h \mbox{ in } \mfa,  \mbox{ with multiplicity } r-1.\\
&&(ii) -\la^2(h),  \mfg_{\la},  \mbox{ with multiplicity } m_{\la},
\mbox{ where } \la\in\Sigma^+. \eeqs
 \par
 Since $|h|=1$,  we have
$\sum_{i}\la^2_i(h)=\fr 12$, which shows that $Ric=-\fr 12$.
\par
We see that at $x=exp_o(r(x)h)$, the eigenvalues with multiplicity
of Hess(r) are
$$0, r-1;\q  |\la(h)|, m_\la.$$
\par
Now we give the rule for the vanishing of  harmonic 1-form in the cases of
irreducible symmetric spaces of noncompact type.
From (\ref{cond2}) we see if $\fr BA\ge 2$, then any harmonic 1-form vanishes, where
$B$ is the Ricc curvature and A is the lower bound of the sectional curvature,
this is proved in \cite{JX}.
\par
Moreover, from (\ref{cond1}), if for any restricted root $\la$, there exists other
two restricted roots $\nu,\mu$   such that
\beq \label{cond3}
|\la|\le |\nu|+|\mu|.
\eeq
Then also any harmonic 1-form vanishes.
\par
We note that in condition (\ref{cond3}) we count the restricted roots with multiplicities.
\par

We are ready to prove Theorem\ref{thmm}.
By Theorem\ref{thmc}, we only consider one of the following  cases,\\
$
SL(3,\mbbr)/SO(3), \qq SU(1,2)/S(U(1)\tm U(2)), \\
SO_0(2,3)/SO(2)\tm SO(3), \qq Sp(2,\mbbr)/U(2).
$\par
We adopt the convention from \cite{HEL}.
Now we verify  the condition (\ref{cond1}) for these cases.
\par
(1)$M=SL(3,\mbbr)/SO(3)$, then $rank(M)=2, dim M=5$,
$\mfg^\mbbc=\mfa_2$, the restricted root system is also $\mfa_2$.
The positive restricted roots are
$$\la_1, \la_2, \la_1+\la_2.$$
We have
$$
|\la_1|\le |\la_2|+|\la_1+\la_2|, |\la_2|\le |\la_1|+|\la_1+\la_2|, |\la_1+\la_2|\le |\la_1|+|\la_2|.
$$
Then condition (\ref{cond1})  is true.
\par
(2) $M=SU(1,2)/S(U(1)\tm U(2))$, then $rank(M)=1, \q  dim M=4$,
$\mfg^\mbbc=\mfa_2$, the  positive restricted roots are
$$
\la_1, \la_1, 2\la_1.
$$
We have
$$
|2\la_1|\le |\la_1|+|\la_1|.
$$
Then condition (\ref{cond1})  is true.
\par
(3)$SO_0(2,3)/SO(2)\tm SO(3)$, then $rank(M)=2, dim(M)=6$,
$\mfg^\mbbc=\mfb_2$, the  positive restricted roots are
$$
\la_1, \la_2, \la_1+\la_2, \la_1+2\la_2.
$$
From
\beqs
&&|\la_1|\le |\la_2|+|\la_1+\la_2|, |\la_2|\le |\la_1|+|\la_1+\la_2|\\
&& |\la_1+\la_2|\le |\la_1| +|\la_2|, |\la_1+2\la_2|\le |\la_2|+|\la_1+\la_2|.
\eeqs
We see that  condition (\ref{cond1})  is true.
\par
(4)
$M=Sp(2,\mbbr)/U(2)$, then $rank(M)=2, dim M=6$,
$\mfg^\mbbc=\mfc_2$,
the  positive restricted roots are
$$
\la_1, \la_2, \la_1+\la_2, 2\la_1+\la_2.
$$
From
\beqs
&&|\la_1|\le |\la_2|+|\la_1+\la_2|, |\la_2|\le |\la_1|+|\la_1+\la_2|\\
&& |\la_1+\la_2|\le |\la_1| +|\la_2|, |2\la_1+\la_2|\le |\la_1|+|\la_1+\la_2|.
\eeqs
We see that  condition (\ref{cond1})  is true.
\par
Now we get Theorem \ref{thmm} from Theorem \ref{thma} and Theorem \ref{vanish2}.

\par
\brem
The only exceptional case is
$SO_0(2,2)/SO(2)\tm SO(2)$, then $rank(M)=2, dim(M)=4$,
$\mfg^\mbbc=\mfd_2$, the  positive restricted roots are
$$
\la_1, \la_2,  \la_1\perp \la_2.
$$
The  condition (\ref{cond1})  is not satisfied.
\erem

\end{document}